\documentclass[12pt,a4paper,twoside,final,notitlepage, leqno]{article}
\usepackage[english]{babel}
\usepackage[T1]{fontenc}
\usepackage{epsfig, graphicx, amssymb}
\usepackage{amsmath,amsthm,epsfig,amsfonts}
\usepackage{float}
\usepackage{color}
\def\br {\break}

\newcommand{\moneq}{\vspace*{-7pt} \begin{equation} \displaystyle }
\newcommand{\moneqstar}{\vspace*{-6pt} \begin{equation*} \displaystyle }
\newcommand{\monendstar}{\vspace*{-6pt} \end{equation*}   }
\newcommand{\monend}{\vspace*{-7pt} \end{equation}   }
\newcommand{\moneqarraystar}{ \begin{eqnarray*} \displaystyle }
\newcommand{\monendarraystar}{ \end{eqnarray*}   }

\newcommand{\dd}{{\rm d}}

\newcommand{\RR}[0]{\mathbb{R}}

\definecolor{vertfonce}{rgb}{0.0, 0.5, 0.0}

\def\section*#1{}

\usepackage{fancyhdr}
\renewcommand{\headrulewidth}{0pt}

\begin{document}

\fancypagestyle{plain}{ \fancyfoot{} \renewcommand{\footrulewidth}{0pt}}
\fancypagestyle{plain}{ \fancyhead{} \renewcommand{\headrulewidth}{0pt}}

~

  \vskip 2.1 cm

\centerline {\bf \LARGE Concavity of thermostatic entropy}

\bigskip

\centerline {\bf \LARGE and convexity of Lax's mathematical entropy}

 \bigskip  \bigskip \bigskip \bigskip \bigskip

\centerline { \large    Fran\c{c}ois Dubois$^{a}$}

\smallskip \smallskip 

\centerline { \it  \small
  $^a$   Aerospatiale,  d\'epartement de math\'ematiques appliqu\'ees, Les Mureaux, France}


\bigskip  \bigskip

\centerline {12 February 1990$\,$\footnote {\rm  \small  $\,$ 
  Published in {\it La Recherche A\'erospatiale}, ISSN 0379-380X, volume 1990-3, pages 77-80, 1990;
edition december 2023.}}




\bigskip  \bigskip
\noindent {\bf \large Abstract}

\noindent
It is shown that for any real gas in chemical equilibrium, the convexity proposed by Lax is equivalent
to the concavity of the entropy considered as thermodynamic potential
and the positivity of temperature. 

\bigskip  \bigskip
\noindent {\bf \large R\'esum\'e}

\noindent
Nous montrons  que pour tout gaz r\'eel \`a l'\'equilibre chimique, la convexit\'e introduite  par Lax
est \'equivalente à la concavit\'e de l'entropie consid\'er\'ee comme potentiel thermodynamique
et \`a la positivit\'e de la temp\'erature.

\newpage 
\bigskip \bigskip    \noindent {\bf \large    1) \quad  Introduction}   

\noindent
Gas dynamics is based on conservation laws of physics
and an assumption of local thermostatic  equilibrium (for instance, see Germain \cite{Ge73})
leading to the Euler equations in the case where viscosity and heat conduction
are neglected. These equations are included in the general category of nonlinear
hyperbolic systems for which Lax suggested a general concept of mathematical
entropy in \cite{FL71,La71}.
In the case of a polytropic ideal gas, it is well known \cite{Ha83,HFM86,Ta84}
that this concept can be interpreted by the usual thermostatic concept of entropy.
In this note, we show that for any real gas in chemical equilibrium, the convexity
proposed by Lax is equivalent to the concavity of the  entropy considered as a 
thermodynamic potential. Below, we recall the bases of thermostatics and gas dynamics in the
next two sections (in order to define the mathematical framework of our analysis),
then, in the last section, we analyse the relations between thermostatic entropy and
Lax's entropy.

\bigskip   \noindent {\bf \large    2) \quad  Thermostatic equilibrium}  

\noindent
In this section, we recall a few essential properties of real gases in thermochemical equilibrium.
These properties are established in the classical works
of Landau-Lifchitz \cite{LL67}, Callen~\cite{Ca66} and Germain {\cite{Ge73}, for instance. 

\smallskip  \noindent
The thermodynamic properties of a gas with mass $ \, M $, internal energy $ \, E \, $ enclosed
in volume~$ \, V \, $ are completely determined by these three parameters
$ \, (M,\, V, E) $. In effect, the thermostatic entropy $ \, S \, $ is a function of the triplet
$ \, (M,\, V, E) $:
\moneq \label{1}
S = \Sigma  \,(M,\, V, \,E) \, .
\monend
We assume that $ \, \Sigma \, $ can be differentiated; its first derivatives
(and second derivatives when defined) are used to compute all the thermodynamic properties
of the gas.

\smallskip  \noindent
We recall the two fundamental properties of entropy $ \, \Sigma $.

\smallskip  \noindent
{\bf Hypothesis 1}.
Function  $ \, \Sigma \, $ is first-order homogeneous (entropy is an extensive variable):
\moneq \label{2}
\Sigma\,  (\lambda \, M,\, \lambda \, V, \,\lambda \, E) = \lambda \, \Sigma \, (M,\, V, \,E)  ,\,
\,\,\,\, \forall \,  \lambda > 0 \, .
\monend

\smallskip  \noindent
{\bf Hypothesis 2}.
Function  $ \, \Sigma \, $ is superadditive:
\moneq \label{3}
 \Sigma \, (M+M',\, V+V', \, E+E') \geq \Sigma \, (M,\, V, \,E) + \Sigma \, (M',\, V', \,E') \, . 
\monend
This last property expresses the second principle of thermodynamics: when two masses of gas
are mixed, the entropy of the resulting system is always higher than or equal to the sum of the
entropies of the constituents. From this, we easily infer the following result.

\smallskip  \noindent
{\bf Proposition 1}.
Function  $ \,\, ]0,\, +\infty [^3 \, \ni  (M,\, V, \,E)
\,\longmapsto \, \Sigma \, (M,\, V, \,E) \,\, $ is concave. 

\smallskip  \noindent
We remark that temperature $ \, T \, $ is the inverse of the derivative of the entropy with respect
to the internal energy:
\moneq \label{4}
{{1}\over{T}} = {{\partial \Sigma}\over{\partial E}}\, (M,\, V, \,E) \, . 
\monend
%

\fancyhead[EC]{\sc{Fran\c{c}ois Dubois}}
\fancyhead[OC]{\sc{Concavity of thermostatic entropy and convexity of Lax's  entropy}}
\fancyfoot[C]{\oldstylenums{\thepage}}

\newpage 
\bigskip    \noindent {\bf \large    3) \quad  Gas dynamics}  

\noindent
The Euler gas dynamics equations express conservation of mass, momentum and energy. They take the
form of a nonlinear hyperbolic system of conservation laws:
\moneq \label{5}
{{\partial U}\over{\partial t}} + {\rm div} \, f(U) = 0 \, .  
\monend
Below, in order to simplify the notations, we will restrict ourselves to the case of a single
space variable $ \, x $, but the extension to more than one dimension is straightforward.
The conserved variables $ \, U \, $ and the flux $ \, f(U) \, $ are expressed as follows
(see for instance Germain~\cite{Ge73}):
\moneq \label{6}
U = \Big( \rho, \, q\equiv \rho \, u ,\,
\varepsilon \equiv \rho \, e + {1\over2} \, \rho \, u^2 \Big)^{\rm t} , 
\monend
\vskip -.5 cm
\moneq \label{7}
f(U) = \Big( \rho \, u ,\, \rho \, u^2 + p,\, \varepsilon \, u + p\, u  \Big)^{\rm t} \,. 
\monend
The pressure $ \, p\, $ introduced in the expression of the flux (\ref{7}) is defined
by the relation\br
$  p = (\gamma-1) \, \rho \, e \,\, $
in the case of an ideal gas with constant specific heats in a ratio $ \, \gamma \, $
(polytropic gas). For a real gas, we have the following general result. 

\smallskip  \noindent
{\bf Proposition 2}.
The pressure $ \, p \, $ and the specific entropy $ \, s \equiv {{\Sigma}\over{M}} \, $             
are defined as functions of the conservative variables $ \, U \, $ alone.

\smallskip  \noindent
{Proof of Proposition 2}.
The second identity of (\ref{6}) allows the internal specific energy $ \, e \, $ to be defined as
a function of $ \, U $:
\moneq \label{8}
e = {{\varepsilon}\over{\rho}} - {{q^2}\over{2 \, \rho}} \, .  
\monend
In addition, considering the homogeneity of $ \, \Sigma $, the specific entropy $ \, s \, $
depends only on the density and the specific internal energy:
\moneq \label{9}
s = \sigma(\rho,\, e) \equiv \Sigma \Big( 1,\, {{1}\over{\rho}} ,\, e \Big)
= {{1}\over{M}} \, \Sigma \, \Big( M,\, {{M}\over{\rho}} ,\, M\, e \Big) \, .
\monend
We then define the pressure by the usual thermostatic properties:
\moneq \label{10}
p = T \, {{\partial \Sigma}\over{\partial V}}  \Big(  1,\, {{1}\over{\rho}} ,\, e \Big) \, . 
\monend
This last expression can be evaluated by function $ \, \sigma \, $ alone and the
conservative variables by the following relation:
\moneq \label{11}
p = - \rho^2 \,\,  {{\displaystyle{{\partial \sigma}\over{\partial \rho}} (\rho,\, e)}\over
{\displaystyle{{\partial \sigma}\over{\partial e}} (\rho,\, e)}} \, . 
\monend
\hfill $\square$

\newpage 
\noindent {\bf \large    4) \quad  Mathematical entropy}  

\smallskip \noindent
In \cite{FL71,La71}, P. Lax proposes a concept of mathematical entropy for any
hyperbolic system of conservation laws (\ref{5}).
It is a function $ \, \eta(U) \, $ which has the following properties:
there exists an entropy flux $ \, \xi(U) \, $ such that
\moneq \label{12}
\dd \xi(U) = \dd \eta(U) \,.\, \dd f(U) \,\,\,\, {\textrm {for any}} \,\, U 
\monend
\vskip -.5 cm 
\moneq \label{13}
U \, \longmapsto \, \eta(U) \,\,  \,\,  {\textrm {is convex}}  . 
\monend
For any hyperbolic system of conservation laws (\ref{1}), nothing {\it a priori}
ensures the existence or uniqueness of a mathematical entropy $ \, \eta \, $
satisfying relations (\ref{12})  and (\ref{13}).

\smallskip \noindent
In the case of gas dynamics, a function $ \, \eta(U) \, $ satisfying equations
(\ref{12}), (\ref{13}) and which is nonaffine has the form
\moneq \label{14}
\eta(U) = -\rho \, s(U) 
\monend
as is suggested by Friedrichs and Lax in \cite{FL71},
where $ \, s(U) \, $ was defined in (\ref{8}) and (\ref{9}).
Below, we shall call this ``Lax's entropy''. 

\smallskip \noindent
Property (\ref{12}) results from the additional conservation law
\moneq \label{15}
{{\partial}\over{\partial t}} \big(\rho \, s(U) \big)
+  {{\partial}\over{\partial x}} \big(\rho \, u \, s(U) \big) = 0 
\monend
satisfied for regular solutions of (\ref{5}), as remarked by Godunov in  \cite{Go59}.

\smallskip \noindent
The convexity of Lax's entropy (\ref{14}) has been thoroughly studied
(for instance Harten \cite{Ha83}, Hughes, Franca and Mallet \cite{HFM86} and
Tadmor \cite{Ta84}) in the case of a polytropic ideal gas for which function $ \, \Sigma \, $
has the following particular form
\moneq \label{16}
\Sigma \,  (M,\, V, \,E) = M \, C_v \, \bigg( 
\log {{E \, M_0}\over{E_0 \, M}} + (\gamma-1)\, \log {{V \, M_0}\over{V_0 \, M}} \bigg) \, . 
\monend
In the case of any real gas, Wagner demonstrated  \cite{Wa87} that this property
was equivalent to the property of convexity of the opposite of the specific entropy
($-s$) with respect to the variables
$\,\, \big( {{1}\over{\rho}} ,\, u ,\, e + {{u^2}\over{2}} \big) \,\, $ 
of the Lagrangian gas dynamics.
But, at our knowledge, the convexity of Lax's entropy does not seem to have been established
generally for suitable physical hypotheses.
When thermochemical equilibrium is assumed reached, we have the following property.

\smallskip  \noindent
{\bf Proposition 3}.
The convexity of function $ \, \eta \, (U) \, $ defined in (\ref{14})
is equivalent to the concavity of function $ \, \Sigma \,  (M,\, V, \,E) \, $ 
and the positivity of the temperature. 

\smallskip  \noindent
{Proof of Proposition 3}.
We first assume $ \, \Sigma \, $ to be concave and $ \, T \, $ to be positive,
{\it i.e.} $ \,  \Sigma  (M,\, V, \, .) \, $ to be nondecreasing. We compute Lax's entropy  (\ref{14})
taking  (\ref{9}) then (\ref{2}) into account. This yields:
\moneqstar 
\eta(U) = -\rho \, s(U) =   -\rho \, \Sigma \, \Big( 1,\, {{1}\over{\rho}} ,\, e \Big) 
\monendstar 
\moneq \label{17}
\eta(U) = - \Sigma \, (\rho,\, 1 ,\, \rho \, e ) \,.
\monend
We define two states $ \, U_1 ,\ U_2 $, and we study
$ \, \eta( (1-t)\, U_1 + t \, U_2 ) $, $\, 0 \leq t \leq 1 $.
With the notations introduced in  (\ref{6}), the internal energy $ \, \rho \, e(t) \, $
associated with state $ \, (1-t)\, U_1 + t \, U_2 \, $ equals, considering  (\ref{8}): 
\moneq \label{18}
\rho \, e(t) = (1-t)\, \varepsilon_1 + t\, \varepsilon_2 - {1\over2} \,
{{((1-t)\, q_1 + t \, q_2)^2}\over{(1-t)\, \rho_1 + t \, \rho_2}} 
\monend
and, remarking that mapping $ \,\, \, ]0,\, +\infty[ \times \RR \ni
  (\rho,\, q) \, \longmapsto \, {{q^2}\over{\rho}} \,\,\,  $ is convex, we deduce
\moneq \label{19}
\rho \, e(t) \geq  (1-t)\, \varepsilon_1 + t\, \varepsilon_2 - {1\over2}
\, \Big( (1-t) \, {{q_1^2}\over{\rho_1}} + t \,  {{q_2^2}\over{\rho_2}} \Big) \, . 
\monend
As function $ \, \Sigma \, $ is non decreasing with respect to the third argument, we deduce
the following inequality from (\ref{19}) 
\moneq \label{20}
\eta \, \big( (1-t)\, U_1 + t \, U_2 \big) \leq
-\Sigma \, \bigg( (1-t)\, \rho_1 + t \, \rho_2 ,\, 1 ,\,
(1-t) \, \Big(\varepsilon_1 - {{q_1^2}\over{2 \, \rho_1}} \Big) +
t \, \Big(\varepsilon_2 - {{q_2^2}\over{2 \, \rho_2}} \Big) \bigg) 
\monend
which, considering the concavity of $ \, \Sigma \, $ and relation  (\ref{17}),
establishes the concavity of Lax's entropy.

\smallskip \noindent
We now assume $ \, \eta \, $ to be convex. We first demonstrate the concavity of $ \, \Sigma $. 
We choose two states $ \, U_1 \, $ and $ \, U_2 $,  with zero momentum:
\moneq \label{21}
q_1 = q_2 = 0 \, . 
\monend
Inequalities (\ref{19}) and (\ref{20}) are then equalities and we deduce the following estimate from
the convexity  of $ \, \eta $:
\moneq \label{22}
\Sigma \, \big(  (1-t)\, \rho_1 + t \, \rho_2 ,\, 1 ,\, (1-t)\, \varepsilon_1 + t \, \varepsilon_2 \big)
\geq (1-t)\, \Sigma \, (\rho_1   ,\, 1 ,\, \varepsilon_1)
+ t \,  \Sigma \, (\rho_2   ,\, 1 ,\, \varepsilon_2) 
\monend
which, considering the homogeneity of $ \, \Sigma $, establishes the required inequality
through an elementary calculation detailed below:

\smallskip \noindent $
\Sigma \big(  (1-t)\, M_1 + t \, M_2 ,\, (1-t)\, V_1 + t \, V_2 ,\, (1-t)\, E_1 + t \, E_2 \big) $

\smallskip \noindent \qquad $
= \big( (1-t)\, V_1 + t \, V_2 \big) \,
\Sigma \, \Big( {{(1-t)\, V_1}\over{ (1-t)\, V_1 + t \, V_2}} \, {{M_1}\over{V_1}}
+ {{t \, V_2}\over{ (1-t)\, V_1 + t \, V_2}}  \, {{M_2}\over{V_2}} ,\, 1 ,\,
 {{(1-t)\, V_1}\over{ (1-t)\, V_1 + t \, V_2}} \, {{E_1}\over{V_1}}
+ {{t \, V_2}\over{ (1-t)\, V_1 + t \, V_2}}  \, {{E_2}\over{V_2}} \Big)  $

\smallskip \noindent \qquad $
\geq  \big( (1-t)\, V_1 + t \, V_2 \big) \, \Big[ 
{{(1-t)\, V_1}\over{ (1-t)\, V_1 + t \, V_2}} \, \Sigma \,\big( {{M_1}\over{V_1}} ,\, 1 ,\, {{E_1}\over{V_1}} \big)
+ {{t \, V_2}\over{ (1-t)\, V_1 + t \, V_2}} \, \Sigma \,\big( {{M_2}\over{V_2}} ,\, 1 ,\, {{E_2}\over{V_2}} \big) \Big] $

\smallskip \noindent \qquad $
= (1-t) \, \Sigma \, (M_1,\, V_1,\, E_1) + t \, \Sigma \, (M_2,\, V_2,\, E_2) \,.$ 

\smallskip \noindent 
We now show that the temperature is positive, {\it i.e.} that $ \, \Sigma \, (M,\, V,\, .) \, $
is a nondecreasing function.
For this, we introduce $ \, \Delta E > 0 \, $ and we define two states, $ \, U_1 \, $
and $ \, U_2 $, as follows:
\moneq \label{23} \left\{ \begin{array} {c}
\rho_1 = \rho ,\, q_1 = \sqrt{8 \, \rho \, \Delta E} ,\, \varepsilon_1 = E + 4 \, \Delta E \\ 
\rho_2 = \rho ,\, q_2 = 0  ,\, \varepsilon_2 = E  \, .
 \end{array} \right.  \monend 
State $ \, (1-t)\, U_1 + t \, U_2 \, $ has an internal energy per unit volume computed according to
relation~(\ref{18}), {\it i.e.}
\moneq \label{24}
\rho \, e(t) = E + 4 \, t \, (1-t) \,  \Delta E \, . 
\monend 
Futhermore, considering   (\ref{17}) and  (\ref{23}), we have 
\moneq \label{25}
(1-t) \, \eta \, (U_1) + t \, \eta \,  (U_2) =  -\Sigma \, (\rho ,\, 1 ,\, E ) , \,\,\,\,
\forall \, t \in [0,\, 1] \, . 
\monend 
The convexity of $ \, \eta \, $ expressed by taking $ \, t = {1\over2} \, $
therefore results in the inequality
\moneq \label{26}
\Sigma \, (\rho ,\, 1 ,\, E + \Delta E) \geq \Sigma \, (\rho ,\, 1 ,\, E)
\monend 
which, considering the homogeneity  (\ref{2}), establishes the required property.
\hfill $\square$

\smallskip \smallskip  \noindent
{\bf Remark}. 
In \cite{Du88}, we demonstrated the following less general result:
if the temperature is positive and the specific internal energy $ \, e \, $
is a regular function (of class $ \, {\cal C}^2 $, excluding any phase transition)
of the (specific entropy $\, s $, specific volume $ \, \tau \equiv {{1}\over{\rho}}$)
pair, the convexity of function
$ \, (s,\, \tau) \, \longmapsto \, e(s,\, \tau) \, $ is equivalent
to that of Lax's entropy. The property of positivity of the temperature
for gases in chemical equilibrium is a classical result.
For the establishment of this property based on purely physical arguments, the reader
is referred, for instance, to Landau and Lifchitz \cite{LL67}.
This last property was also used by Friedrichs and Lax in \cite{FL71}.
The convexity of Lax's entropy was probably already known by Godunov, who introduced in 
\cite{Go59} the entropy variables $ \, \varphi $, {\it i.e.} the gradient of  (\ref{14})
with respect to the conservative variables (the bijectivity of the mapping
$ \, U \,\longmapsto \, \varphi \, $ results immediately from the strict
convexity of $ \, \eta $).
For an extension of the analysis concerning mixing of two gases without mutual
interaction, we refer the reader to  \cite{Du89}.

\bigskip \bigskip      \noindent {\bf  \large  References }

\fancyfoot[C]{\oldstylenums{\thepage}}



\end{document}